\documentclass[a4paper,notitlepage,twoside,reqno,11pt]{amsart}

\usepackage{amsmath,amssymb,amscd,amsthm,amsfonts}
\usepackage{mathrsfs}
\usepackage{latexsym}
\usepackage{bbm,pifont}

\usepackage{anysize}
\marginsize{3.4cm}{3.4cm}{3cm}{3cm}

\usepackage{color,xcolor,graphicx,graphics}
\usepackage{subfigure}
\usepackage{tikz,tikz-cd}
\usepackage{extarrows}

\usepackage{dcolumn}
\usepackage{indentfirst}
\usepackage{enumerate,enumitem}

\usepackage{titletoc}
\usepackage[colorlinks=false,linkcolor=black,citecolor=black,urlcolor=black,bookmarksopen=true,linktocpage=true,pdfstartview=FitH]{hyperref}
\usepackage{caption2}

\usepackage{xy}

\newtheorem{theorem}{Theorem}[section]
\newtheorem{lemma}[theorem]{Lemma}

\theoremstyle{remark}
\newtheorem{remark}[theorem]{Remark}

\newcommand{\C}{\mathbb C}
\newcommand{\Chat}{\widehat{\mathbb C}}
\newcommand{\D}{\mathbb D}
\newcommand{\T}{\mathbb T}

\makeatletter
\@addtoreset{equation}{section}
\makeatother

\begin{document}

\title{Quasicircularity of the Zakeri Curve and Capture Component Boundaries}

\author{YUMING FU}
\address{School of Mathematical Sciences, Shenzhen University, Shenzhen, 518060, P. R. China}
\email{yumingfuxy@szu.edu.cn}

\author{WEIWEI SU}
\address{School of Mathematical Sciences, Nanjing University, Nanjing, 210093, P. R. China}
\email{suweiwei2024@163.com}

\begin{abstract}
For a fixed irrational rotation number of bounded type, we prove that the locus of critically marked cubic polynomials whose two critical points lie on the boundary of the fixed Siegel disk is a quasicircle. This locus is known as the Zakeri curve. The main step is a boundary correspondence lemma for holomorphic families of disks with uniformly quasiconformal extensions. We further prove that every capture component in this family is a quasidisk, including components attached to the Zakeri curve.
\end{abstract}

\date{\today}

\maketitle

\section{Introduction}\label{sec:intro}

Parameter loci defined by critical points on the boundary of a linearization domain describe how the critical points participate in the dynamics as the parameter varies. A classical example arises in cubic polynomials with a fixed bounded-type Siegel disk: both critical points lie on the boundary of the Siegel disk, and the parameters with this property form a curve \cite{Zak99}. Its topological and geometric properties are natural questions in the study of the corresponding parameter slice.

Fix an irrational number $\theta$ of bounded type, meaning that the coefficients of its continued fraction expansion are bounded, and set $\lambda=e^{2\pi i\theta}$. We use the critically marked normalization
\begin{equation*}\label{eq:family}
 P_c(z)=\lambda z\left(1-\frac{1+c^{-1}}2z+\frac{z^2}{3c}\right),
\end{equation*}
where $c\in\C^*$. Its critical points are $1$ and $c$. Let $\Delta_c$ denote the Siegel disk centered at $0$, and set
\[
 \Gamma=\{c\in\C^*:1,c\in\partial\Delta_c\}.
\]
Zakeri proved that $\Gamma$ is a Jordan curve separating $0$ from $\infty$ \cite[Theorem~14.3]{Zak99}. We call it the \emph{Zakeri curve}. His surgery construction also gives uniform quasidisk geometry for the Siegel disks in the dynamical planes. The question considered here is whether the parameter curve $\Gamma$ is a quasicircle.

A question recorded at the 2008 Oberwolfach meeting asks whether the Zakeri arc is a quasiarc when the rotation number is of bounded type \cite[p.~2708, Question~(13)]{PLSS08}. Interchanging the two critical points of $P_c$ identifies $c$ with $1/c$; the quotient coordinate is $(c+c^{-1})/2$. Thus the arc formulation in unmarked moduli space and the present curve formulation use different critical markings. Our theorem also implies that this quotient arc is a quasiarc, as explained in Remark \ref{rem:unmarked}, and gives an affirmative answer to this question.

Related Jordan-curve results occur in other Siegel parameter spaces. For quadratic rational maps with a fixed Siegel disk with bounded-type rotation number, Zhang proved that the locus where both critical points lie on the boundary of the Siegel disk is a Jordan curve \cite{Zha08b}. For quadratic rational maps with a two-cycle of Siegel disks with a bounded-type rotation number, Fu, Yang and Zhang proved that the locus where each boundary of the two-cycle of Siegel disks contains exactly one critical point is a Jordan curve \cite{FYZ22}. For cubic polynomials with a two-cycle of Siegel disks, Fu, Hu and Muzician proved that the locus with one critical point on each of the two boundaries is a Jordan curve. There are also two Jordan arcs corresponding to both critical points lying on the boundary of the same disk \cite{FHM24}. These results show that simultaneous critical boundary conditions lead to natural curves and arcs in several one-dimensional parameter spaces.

An attracting analogue requires specifying the domain whose boundary is used. Building on the analytic model of Petersen and Tan \cite{PT09}, Ch\'eritat proved that, in a cubic slice with fixed multiplier $0<|\lambda|<1$, the parameters for which both critical points lie on the boundary of the maximal Koenigs linearization domain form a Jordan curve \cite{Che20B}. Here the domain is the image of the disk of convergence of the inverse of the linearizing map. It is contained strictly inside the attracting basin, so its boundary is different from the boundary of the immediate basin. For an example involving an actual attracting basin boundary, Qiu, Roesch, Wang and Yin proved that the outer escape component of the McMullen family $z\mapsto z^n+\kappa/z^n$, $n\geq3$, has Jordan-curve boundary \cite{QRWY15}. Their boundary characterization allows either the free critical points or a parabolic cycle on the boundary of the immediate basin of infinity.

Quasicircle results in parameter space precede the present work. A capture component is a connected component of the set of parameters for which one critical orbit enters $\Delta_c$. In the same cubic Siegel slice, Wang and Yang proved that capture components have Jordan boundaries and that a capture component is a quasidisk when its boundary is disjoint from $\Gamma$ \cite{WY23}. The boundary of each capture component is locally a quasiarc away from $\Gamma$. The holomorphic motion used in their argument is available on either side of $\Gamma$, but that argument does not directly control attachment points on $\Gamma$. We first prove that $\Gamma$ is a quasicircle with a quasisymmetric critical-angle parametrization. We then use this geometry to treat the attachment point of a capture component and obtain quasicircularity of every capture boundary in the critically marked family.

Write $I$ and $E$ for the components of $\Chat\setminus\Gamma$ containing $0$ and $\infty$, respectively. For $c\in\overline E\setminus\{\infty\}$, normalize the inverse of the linearizing map by
\begin{equation*}\label{eq:linearizer}
 h_c:\D\longrightarrow\Delta_c, \quad \text{with} \quad h_c(0)=0,
\end{equation*}
 and an extension 
$\hat{h}_c:\overline{\D}\to\overline{\Delta_c}$
satisfying
\begin{equation*}
   \hat{h}_c(1)=1 \quad \text{and} \quad  P_c\circ \hat{h}_c=\hat{h}_c\circ R_\theta,
\end{equation*}
where $R_\theta(z)=\lambda z$. The boundary extension exists because $\Delta_c$ is a quasidisk. Let $\T=\partial\D$. The critical-angle coordinate is
\begin{equation*}\label{eq:alpha}
 \alpha:\Gamma\longrightarrow\T,\qquad\text{and}\qquad \alpha(c)=\hat{h}_c^{-1}(c).
\end{equation*}

\begin{theorem}\label{thm:main}
The curve $\Gamma$ is a $K(\theta)$-quasicircle. Moreover, if $g:\Chat\setminus\overline\D\to E$ is a conformal map with an extension $\hat{g}:\Chat\setminus\D \to \overline{E}$, then
\[
 \alpha\circ \hat{g}|_{\T}:\T\longrightarrow\T
\]
is quasisymmetric.
\end{theorem}

\begin{theorem}\label{thm:capture}
Every capture component of the critically marked family $P_c$ defined above is a quasidisk. Equivalently, its boundary is a quasicircle.
\end{theorem}

The proof is based on Lemma \ref{lem:mainlemma} in Section~\ref{sec:one}. It concerns a holomorphic family of disk embeddings uniform quasiconformal extensions of the maps in this family, and a boundary homeomorphism between the disk boundaries and the parameter boundary. The argument uses a cross ratio and harmonic measure to obtain a lower bound on the lengths of arcs with fixed harmonic mass. Varying the pole gives control of every balanced quadruple and hence yields quasisymmetry. Section~3 applies this lemma using the uniform quasicircle property of the Siegel disks and Zakeri's continuity results, together with the topological parametrization of Yang and Zhang \cite{YZ24}. The inversion symmetry of the cubic family then turns the boundary quasisymmetric map into a quasiconformal map of the sphere. Section~\ref{sec:capture} proves Theorem~\ref{thm:capture}. We detail the degree-one argument at the level-one root $c=1$. A zero-counting argument then proves uniqueness of the exterior capture component at each higher-level root, allowing the same degree-one proof there.

\section{Quasisymmetric boundary correspondences}\label{sec:one}

\subsection{Quasiconformal preliminaries}

We state the four classical results used below, including a cross-ratio distortion theorem for quasiconformal maps. For four distinct points $a,b,c,d$ of $\Chat$, denote their cross-ratio by
\[
 [a,b;c,d]=\frac{(a-b)(c-d)}{(a-c)(b-d)},
\]
with $[\infty,b;c,d]=\frac{c-d}{b-d}$.

\begin{theorem}[Cross-ratio distortion]\label{thm:distortion}
For every $K\geq1$, there is an increasing homeomorphism $\eta_K:[0,\infty)\to[0,\infty)$ such that every $K$-quasiconformal homeomorphism $H$ of $\Chat$ satisfies
\begin{equation}\label{eq:distortion}
 |[H(a),H(b);H(c),H(d)]|
 \leq \eta_K\bigl(|[a,b;c,d]|\bigr)
\end{equation}
for all distinct points $a,b,c,d\in\Chat$.
\end{theorem}

Equation~\eqref{eq:distortion} follows from Ahlfors \cite[Chapter III, Section C, Theorem 2]{Ahl06}. For the $K$-quasiconformal maps $G:\C\to\C$ fixing $0$ and $1$, uniform H\"older estimates give $m_K(t):=\sup_G\sup_{|z|\le t}|G(z)|<\infty$ and $m_K(t)\to0$ as $t\to0$. Choose an increasing homeomorphism $\eta_K\ge m_K$. M\"obius normalizations sending $(a,c,d)$ and their images to $(0,1,\infty)$ produce such a $G$ with $G([a,b;c,d])=[H(a),H(b);H(c),H(d)]$. This implies \eqref{eq:distortion}. 

Interchanging the second and third entries replaces the cross ratio by its reciprocal, so
\[
 |[H(a),H(b);H(c),H(d)]|\ge
 \frac{1}{\eta_K(1/|[a,b;c,d]|)}.
\]

We can assume the additional condition $\eta_K(t)>t$ for every $t>0$ by replacing $\eta_K(t)$ with $\eta_K(t)+t$ if necessary.

\begin{theorem}[Beurling--Ahlfors criterion]\label{thm:qs}
Let $\sigma:\T\to\T$ be an orientation-preserving homeomorphism. The following conditions are equivalent:
\begin{enumerate}[label=\textup{(\roman*)}]
\item $\sigma$ is quasisymmetric.
\item $\sigma$ extends to a quasiconformal homeomorphism of $\D$ onto itself.
\item There is $b\in(0,1/2)$ such that
\[
 b\leq[\sigma(x_1),\sigma(x_2);\sigma(x_3),\sigma(x_4)]\leq1-b
\]
whenever $x_1,x_2,x_3,x_4\in\T$ are cyclically ordered and
$[x_1,x_2;x_3,x_4]=1/2$.
\end{enumerate}
\end{theorem}

This is the circle criterion stated in Sergeev \cite[Section 2.2, Theorem~1]{Ser09}. His cross ratio is $1-[x_1,x_2;x_3,x_4]$, and his constant $\varepsilon$ corresponds to $b=(1-\varepsilon)/2$, so the criterion has exactly the form above with our convention. The real-line version and extension theorems are given in Ahlfors \cite[Chapter~IV, Theorems~1 and~2]{Ahl06}. For an orientation-reversing homeomorphism, composition with a reflection gives the same cross-ratio criterion for quasisymmetry.

\begin{theorem}[Sewing across an analytic arc]\label{thm:sewing}
Let $U$ be a plane domain, let $\gamma\subset U$ be an analytic arc, and let $F$ be an orientation-preserving homeomorphism on $U$. If $F$ is $K$-quasiconformal on $U\setminus\gamma$, then it is $K$-quasiconformal on $U$.
In particular, a sphere homeomorphism that is $K_1$-quasiconformal on $\D$ and $K_2$-quasiconformal on $\Chat\setminus\overline\D$ is $\max\{K_1,K_2\}$-quasiconformal on $\Chat$.
\end{theorem}

The analytic-arc sewing theorem is cited as Lehto and Virtanen \cite[Chapter~I, Theorem~8.3]{LV73}. The circle case follows by applying the arc statement to $\T$. This is used in the proof of Theorem \ref{thm:main}.

\begin{theorem}[Bounded-turning criterion]\label{thm:bounded-turning}
Let $\gamma\subset\C$ be a Jordan curve. The following conditions are equivalent:
\begin{enumerate}[label=(\roman*),leftmargin=*,itemsep=4pt]
 \item $\gamma$ is a quasicircle, that is, $\gamma=H(\T)$ for a quasiconformal homeomorphism $H:\C\to\C$.
 \item There is a constant $C\ge1$ such that, for every pair of distinct points $x,y\in\gamma$, the two closed subarcs $\gamma_1,\gamma_2$ of $\gamma$ with endpoints $x,y$ satisfy
 \[
 \min\{\operatorname{diam}\gamma_1,\operatorname{diam}\gamma_2\}\le C|x-y|.
 \]
\end{enumerate}
\end{theorem}

This is the bounded-turning formulation of Ahlfors' characterization \cite[Theorem 1]{Ahl63}. Here diameter is Euclidean, and condition~\textup{(ii)} is called the bounded-turning condition. 

\subsection{A boundary correspondence lemma}

\begin{lemma}\label{lem:mainlemma}
Let $\Omega\subset\Chat$ be a Jordan domain. Suppose that
\[
 h:\overline\Omega\times\overline\D\longrightarrow\Chat
 \qquad\text{and}\qquad
 a:\overline\Omega\longrightarrow\Chat
\]
are continuous and satisfy the following conditions:
\begin{enumerate}
\item For each $p\in\overline\Omega$, the map $h_p=h(p,\cdot)$ is an embedding of $\overline\D$ and admits a $K$-quasiconformal extension to $\Chat$, with a uniform constant $K$ for all $p$.
\item The maps $p\mapsto h_p(z)$, for $z\in\D$, and $p\mapsto a(p)$ are holomorphic on $\Omega$ as maps into $\Chat$.
\item For $p\in\Omega$, one has $a(p)\notin h_p(\overline\D)$. For $p\in\partial\Omega$, one has
\[
 a(p)=h_p(\alpha(p)),
\]
where $\alpha:\partial\Omega\to\T$ is a homeomorphism.
\end{enumerate}
Then for every conformal map $g:\D\to\Omega$ with an extension $\hat{g}:\overline{\D}\to \overline{\Omega}$, $\alpha\circ \hat{g}|_{\T}$ is quasisymmetric. 
\end{lemma}

\begin{proof}
Write $\eta=\eta_K$, with the choice made in Theorem \ref{thm:distortion}.

\medskip\noindent\textbf{Step 1. Normalize at the pole and estimate a cross ratio.}

Fix $p_0\in\Omega$ and choose a $K$-quasiconformal extension $H_{p_0}$ of $h_{p_0}$. Since
\[
 z_0=H_{p_0}^{-1}(a(p_0))\in\Chat\setminus\overline\D,
\]
there is a Möbius transformation $M$ that preserves the unit disk and satisfies $M(\infty)=z_0$. Fix the map $M$ when $p$ varies.

For $u\in\D$ with $3/4\leq|u|<1$, define
\begin{equation*}\label{eq:Fu}
 F_u(p)=[a(p),h_p(M(u));h_p(M(0)),h_p(M(1/2))].
\end{equation*}
The four points are distinct throughout $\overline\Omega$, since the three interior points are distinct, and the fourth point is outside the closed disk for $p\in\Omega$ and on its boundary for $p\in\partial\Omega$. Since the cross ratio depends holomorphically on four distinct points, $F_u$ is a finite, nonvanishing holomorphic function on $\Omega$, continuous and nonvanishing on $\overline\Omega$ for fixed $u$, even if some point is infinite. Consequently $\log|F_u|$ is harmonic, and
\begin{equation}\label{eq:mean}
 \log|F_u(p_0)|=\int_{\partial\Omega}\log|F_u(p)|\,d\omega_{p_0}(p),
\end{equation}
where $\omega_{p_0}$ denotes harmonic measure in $\Omega$ with pole $p_0$.

At the pole, $a(p_0)=H_{p_0}(M(\infty))$ and
\[
 |[\infty,u;0,1/2]|=\frac{|1/2|}{|u-1/2|}\in[1/3,2].
\]
The reciprocal form of the equation \eqref{eq:distortion}, applied to $H_{p_0}\circ M$, gives
\begin{equation}\label{eq:pole}
 \log|F_u(p_0)|\geq-\log\eta(3).
\end{equation}
For $p\in\partial\Omega$, let $t=M^{-1}(\alpha(p))\in\T$. A $K$-quasiconformal extension of $h_p$ gives
\begin{equation}\label{eq:boundary}
 |F_u(p)|\leq\eta\left(\frac{|t-u|\,|1/2|}{|t|\,|u-1/2|}\right)
 \leq\eta(2|t-u|)\leq\eta(4).
\end{equation}
The extension can be chosen independently for each boundary parameter.

\medskip\noindent\textbf{Step 2. Bound below the length of an arc with harmonic mass $1/4$.}

Consider the probability measure
\[
 \mu_{p_0}=(M^{-1}\circ\alpha)_*\omega_{p_0}
\]
on $\T$, where $(\cdot)_*$ denotes pushforward of measures.

We claim that every arc $J\subset\T$ with $\mu_{p_0}(J)=1/4$ has angular length at least $\delta=\delta(K)>0$, independently of the choice of $p_0$.

It suffices to consider $\ell=\operatorname{length}(J)\leq1/4$, measured in radians. If $\xi$ is the midpoint of $J$, take $u=(1-\ell)\xi$. Then $3/4\leq|u|<1$, and for $t=\xi e^{i\beta}\in J$, i.e., $|\beta|\le\ell/2$,
\[
|t-u|\le |t-\xi|+|\xi-u|\le |\beta|+\ell\le 2\ell.
\]
On the corresponding part $A$ of $\partial\Omega$, equation \eqref{eq:boundary} implies $|F_u(p)|\leq\eta(4\ell)$; on the remaining part $B$ of $\partial\Omega$, it gives $|F_u(p)|\leq\eta(4)$. Combining these bounds with equations \eqref{eq:mean} and \eqref{eq:pole} yields
\begin{equation*}\label{eq:quarter}
\begin{aligned}
-\log\eta(3)\le \log|F_u(p_0)|
&=\int_{\partial\Omega}\log|F_u(p)|\,d\omega_{p_0}(p)\\
&=\int_A\log|F_u|\,d\omega_{p_0}
+\int_{B}\log|F_u|\,d\omega_{p_0}\\
&\le\frac14\log\eta(4\ell)
+\frac34\log\eta(4).
\end{aligned}
\end{equation*}
Therefore, 
\[
\ell\ge\frac14\,\eta^{-1}\!\left(\frac{1}{\eta(3)^4\eta(4)^3}\right)>0
\]
since $\eta$ is increasing. This proves the claim by taking $\delta=\frac14\,\eta^{-1}\!\left(\frac{1}{\eta(3)^4\eta(4)^3}\right)$. Decrease $\delta$ if necessary so that $\delta\leq1/4$.

\medskip\noindent\textbf{Step 3. Control balanced quadruples and apply the quasisymmetry criterion.}
Let $\tau=\alpha\circ \hat{g}|_{\T}$, and take four positively cyclically ordered points $x_1,x_2,x_3,x_4\in\T$ such that
\begin{equation*}\label{eq:balanced}
 [x_1,x_2;x_3,x_4]=\frac12.
\end{equation*}
There is a disk automorphism $\phi$ sending $1,i,-1,-i$ to these four points in order: three images determine a circle-preserving M\"obius map, and the cross-ratio determines the fourth. Choose the pole $p_0=g(\phi(0))$. By conformal invariance of harmonic measure \cite[p.~13]{GM05}, each of the four arcs of $\partial\Omega$ between successive $\hat{g}(x_j)$ has $\omega_{p_0}$-mass $1/4$.

Apply Step~2 with this pole. Each of the four successive arcs between
\[
 y_j=M^{-1}(\tau(x_j)), 
\]
for $j=1,2,3,4$ has length at least $\delta$. Every chord between distinct $y_j$ therefore has length at least $2\sin(\delta/2)$. For cyclic points on a circle, their cross ratio $R=[y_1,y_2;y_3,y_4]$ belongs to $(0,1)$, and
\[
 1-R=\frac{(y_1-y_4)(y_2-y_3)}{(y_1-y_3)(y_2-y_4)}.
\]
Using $2\sin(\delta/2)\le |y_i-y_j|\le 2$ for distinct $i,j$ and cross-ratio is invariant under the M\"obius map, we obtain
\begin{equation*}\label{eq:balanced-bound}
 \sin^2(\delta/2)\leq[\tau(x_1),\tau(x_2);\tau(x_3),\tau(x_4)]\leq1-\sin^2(\delta/2).
\end{equation*}
These estimates also hold when $\tau$ reverses cyclic order. Theorem~\ref{thm:qs} now implies that $\tau$ is quasisymmetric.
\end{proof}

\section{Quasicircularity of the Zakeri curve}\label{sec:two}

In this section, we verify that the set $E$ in the parameter space and the map $\hat{h}_c : \overline{\D} \to \overline{\Delta_c}$ defined in Section~\ref{sec:intro} satisfy the hypotheses of Lemma \ref{lem:mainlemma}, including in particular the parameter at infinity. Define
\[
    h_\infty : \mathbb{D} \to \Delta_\infty
\]
to be conformal with an extension $\hat{h}_\infty:\overline{\D}\to \overline{\Delta_\infty}$, where $\Delta_\infty$ is the Siegel disk of the map $Q_\theta(z) = \lambda z(1 - z/2)$ centered at $0$, normalized by
\[
    \hat{h}_\infty(0) = 0, \qquad \hat{h}_\infty(1) = 1,
\]
so that $Q_\theta \circ \hat{h}_\infty(z) = \hat{h}_\infty(\lambda z)$.

\begin{lemma}\label{lem:family}
The maps can be extended to a continuous family on $\overline E\times\overline\D$. For each $z\in\D$, the map $c\mapsto \hat{h}_c(z)$ is holomorphic on $E$. Every $\hat{h}_c$, including $\hat{h}_\infty$, admits a $K(\theta)$-quasiconformal extension to $\Chat$. Furthermore,
$c\notin \hat{h}_c(\overline\D)$ when $c\in E$, while $c\in \hat{h}_c(\T)$ when $c\in \Gamma$.
\end{lemma}

\begin{proof}
Zakeri's surgery in \cite{Zak99} gives a uniform quasiconformal extension bound of $\hat{h}_c$. Meanwhile the result in \cite{Zha11} also implies the boundaries $\partial\Delta_c$ are $K'(\theta)$-quasicircles; the quasidisk extension theorem therefore gives sphere $K(\theta)$-quasiconformal extensions of $\hat{h}_c$, with $K(\theta)$ independent of $c$. 

For finite parameters, Zakeri's boundary continuity theorem \cite[Theorem 13.9]{Zak99}, together with quasiconformal compactness and the normalization of $h_c$, gives continuity of $c\mapsto \hat{h}_c$ in the uniform topology on $\overline\D$.

It remains to treat $c=\infty$. The same theorem in \cite{Zha11} applies to $Q_\theta$, whose unique finite critical point is $1$; thus the boundary of its Siegel disk $\Delta_\infty$ is a quasicircle containing $1$. For any sequence $c_n\to\infty$, choose uniformly quasiconformal extensions $H_n$ of $\hat{h}_{c_n}$. Since
\[
 H_n(0)=0,\quad H_n(1)=1\quad \text{and} \quad
 H_n(\lambda)=P_{c_n}(1)\longrightarrow\lambda/2,
\]
quasiconformal compactness gives a subsequence converging uniformly in the spherical metric to a $K(\theta)$-quasiconformal map $H$. For every fixed integer $m\ge0$,
\[
 H(\lambda^m)=\lim_{n\to\infty}P_{c_n}^{m}(1)
 =Q_\theta^{m}(1)=\hat{h}_\infty(\lambda^m).
\]
Density of $\{\lambda^m:m\ge0\}$ in $\T$ gives $H|_{\T}=\hat{h}_\infty|_{\T}$. Since $H(0)=0$, we have $H(\D)=\Delta_\infty$; moreover, $H|_{\D}$ is conformal as a limit of the maps $\hat{h}_{c_n}$. Uniqueness of the normalized linearization map implies $H|_{\overline\D}=\hat{h}_\infty$. Every convergent subsequence has the same restriction to $\overline{\mathbb D}$, so $\hat{h}_c\to \hat{h}_\infty$ uniformly on $\overline\D$. The map $H$ also supplies the required extension of $\hat{h}_\infty$.

For $c\in E\setminus\{\infty\}$, we have $\hat{h}_c(\lambda^m)=P_c^m(1)$. We know that $P^n_c(1)$ is holomorphic in $c$. For each $\zeta\in \T$, we can take a sequence $\lambda^{n_j}\to \zeta$, such that $P_c^{n_j}(1)=\hat{h}_c(\lambda^{n_j})\to \hat{h}_c(\zeta)$. Since $(c,z)\mapsto \hat{h}_c(z)$ is uniformly continuous on $L\times\T$, where $L$ is any compact set of $E\setminus\{\infty\}$, $\hat{h}_c(\zeta)$ is holomorphic in $c$. For $z\in \D$, $\hat{h}_c$ is holomorphic by Cauchy's formula. Continuity at $\infty$ implies that the singularity there is removable. Thus, $\hat{h}_c$ is holomorphic on $E\times\D$.

Finally, Zakeri's description of the critical points \cite[Corollaries 9.3 and 13.6, Theorem 14.3]{Zak99} shows $c\notin\partial\Delta_c$ on $E\setminus\{\infty\}$; since a critical point cannot lie inside a Siegel disk, $c\notin\overline{\Delta_c}$. On $\Gamma$, the assertion follows from its definition. At $c=\infty$, it follows from boundedness of $\Delta_\infty$.
\end{proof}

\begin{lemma}\label{lem:angle}
The map $\alpha:\Gamma\to\T$ defined in Section~\ref{sec:intro} is a homeomorphism.
\end{lemma}

\begin{proof}
We apply Yang and Zhang \cite[Proposition~5.4 and Theorem~5.7]{YZ24}. Their normalization is
\[
 f_b(z)=\lambda z\left(1-\frac{b+b^{-1}}2z+\frac{z^2}{3}\right),
\]
satisfying $P_{b^2}(bz)=b f_b(z)$. Thus $c=b^2$ identifies their quotient by $b\sim-b$ with our parameter plane; their quadratic model $q_\theta$ is our $Q_\theta$. Theorem~5.7 gives a homeomorphism $\Psi$ whose inverse image $X=\Psi^{-1}(\partial\Delta_\infty)$ is a Jordan curve. For $n\geq0$, put $c_{-n}=\Psi^{-1}(\hat{h}_\infty(\lambda^{-n}))$. Proposition~5.4 gives $c_{-n}\in\Gamma$ and
\[
 P_{c_{-n}}^n(c_{-n})=1,
 \quad \text{and} \quad \alpha(c_{-n})=\lambda^{-n}
       =\hat{h}_\infty^{-1}(\Psi(c_{-n})).
\]
These roots are dense in $X$, so $X\subset\Gamma$, and hence $X=\Gamma$, since both are Jordan curves. Continuity from Lemma \ref{lem:family} now gives $\alpha=\hat{h}_\infty^{-1}\circ\Psi$ on $\Gamma$. Hence $\alpha$ is a homeomorphism.
\end{proof}

\begin{remark}
    The method used in \cite{FYZ22} can also be used to prove that $\alpha$ is a homeomorphism.
\end{remark}

\begin{proof}[Proof of Theorem \ref{thm:main}]
Let $g:\Chat\setminus\overline\D\to E$ be as in the statement, and $\hat{g}$ be its homeomorphic boundary extension. By Lemma \ref{lem:family} and Lemma \ref{lem:angle}, Lemma \ref{lem:mainlemma} can be applied with
\[
 \Omega=E,\quad h(c,z)=\hat{h}_c(z),\quad a(c)=c\quad \text{and} \quad
 \alpha(c)=\hat{h}_c^{-1}(c).
\]
Using $z\mapsto \hat{g}(1/z)$, we obtain the quasisymmetry of $\tau=\alpha\circ \hat{g}|_{\T}$.

To obtain the geometry of $\Gamma$, we use the symmetry of $P_c(z)$:
\begin{equation*}\label{eq:symmetry}
 P_{1/c}(z)=c^{-1}P_c(cz).
\end{equation*}
It follows that $\iota(c)=1/c$ preserves $\Gamma$ and exchanges $I$ and $E$. For $c\in\Gamma$, $|\alpha(c)|=1$, consider the conformal map 
\[
z\mapsto c^{-1}h_c(\alpha(c)z):\D\to\Delta_{1/c}
\]
whose continuous extension fixes $0$ and $1$. Uniqueness of the normalized linearization map gives
\[
 \hat{h}_{1/c}(z)=c^{-1}\hat{h}_c(\alpha(c)z).
\]
Take $z=1/\alpha(c)$, we obtain $\hat{h}_{1/c}(1/\alpha(c))=1/c.$ Since $\alpha(1/c)=\hat{h}_{1/c}^{-1}(1/c)$ and $\hat{h}_{1/c}$ is a homeomorphism on $\T$, we have
\[
\alpha(1/c)=1/\alpha(c).
\]

Assume first that $\tau$ preserves orientation. The exterior-disk form of Theorem~\ref{thm:qs} gives a quasiconformal extension $\widetilde\tau:\Chat\setminus\D\to\Chat\setminus\D$. Define
\begin{equation}\label{eq:directextension}
 \Phi(z)=
 \begin{cases}
 \hat{g}\bigl(\widetilde\tau^{-1}(z)\bigr),&z\in\Chat\setminus\D,\\[2mm]
 \displaystyle\frac{1}{\hat{g}\bigl(\widetilde\tau^{-1}(1/z)\bigr)},&z\in\overline\D.
 \end{cases}
\end{equation}
The boundary values agree, since for $\zeta\in\T$,
\[
 \hat{g}(\tau^{-1}(\zeta))=\alpha^{-1}(\zeta)
 =\frac{1}{\alpha^{-1}(1/\zeta)}
 =\frac{1}{\hat{g}(\tau^{-1}(1/\zeta))}.
\]
The two pieces map $\Chat\setminus\overline\D$ and $\D$ onto $E$ and $I$, respectively. Thus $\Phi$ is a sphere homeomorphism, quasiconformal off $\T$, and Theorem~\ref{thm:sewing} makes it quasiconformal everywhere. Since $\Phi(\T)=\Gamma$, the curve $\Gamma$ is a quasicircle.

If $\tau$ reverses orientation, the argument is the same by replacing $\tau$ by $\bar{\tau}$.
\end{proof}

\begin{remark}[The unmarked locus]\label{rem:unmarked}
Let $q(z)=(z+z^{-1})/2$. In either orientation case, the map $\Phi$ in \eqref{eq:directextension} satisfies $\Phi(1/z)=1/\Phi(z)$. It therefore descends to a sphere homeomorphism $\widehat\Phi$ satisfying $q\circ\Phi=\widehat\Phi\circ q$. This map is quasiconformal off the branch values $\pm1$, and hence is quasiconformal everywhere by Theorem~\ref{thm:sewing}. Since $q(\T)=[-1,1]$, we obtain $q(\Gamma)=\widehat\Phi([-1,1])$, a quasiarc.
\end{remark}

\section{Quasicircularity of capture component boundaries}\label{sec:capture}
Wang and Yang proved that a capture boundary is locally a quasiarc at every point outside $\Gamma$; in particular, a capture component whose boundary is disjoint from $\Gamma$ is already a quasidisk \cite[Theorem 4.2]{WY23}. Together with the symmetry, it therefore suffices to prove the local quasiarc property at the attachment point of a component in $E$.

Throughout this section, fix a level-$n$ capture component $\mathcal C\subset E$ with $n=\min\{k\ge1:P_c^k(c)\in\Delta_c\}$,
whose boundary meets $\Gamma$. It is a Jordan domain, and
\[
 \Phi_{\mathcal C}(c)=h_c^{-1}(P_c^n(c)):\mathcal C\longrightarrow\D
\]
is conformal by \cite[Lemma 2.4 and Theorem 3.4]{WY23}. Its boundary intersects $\Gamma$ at the single point $c_*$, and
\begin{equation*}\label{eq:capture-root}
 P_{c_*}^{n-1}(c_*)=1\quad \text{and} \quad\widehat\Phi_{\mathcal C}(c_*)=\lambda
\end{equation*}
where $\widehat\Phi_{\mathcal C}$ denotes the extension \cite[Lemma 3.2 and Corollary 3.5(a)]{WY23}. Set 
\[
q(c)=P_c^{n-1}(c),\quad \text{for}\  c\in\C^{*},\ \text{and}\  \beta(c)=\lambda^{n-1}\alpha(c),\quad \text{for}\ c\in\Gamma.
\]
Let $\psi(z)=\Phi_{\mathcal C}^{-1}(\lambda z)$ for $z\in\D$. We use the continuous boundary extension of $\psi$, so that $\psi(1)=c_*$. 

On $\Gamma$ one has $q(c)=h_c(\beta(c))$ and $\beta(c_*)=1$. Our only task is to prove that $\partial\mathcal C$ is locally a quasiarc at $c_*$. All estimates below concern this fixed root, and their constants may depend on $\mathcal C$. We first treat $n=1$, then verify the changes needed for $n\ge2$.

\subsection{Harmonic measure near an attachment point}
We use the following local consequence of Lemma \ref{lem:mainlemma} for the map $q(c)$ just defined.

\begin{lemma}\label{lem:local-mass}
Let $\Omega\subset E$ be a small Jordan domain bounded by an arc $A\subset\Gamma$ through $c_*$ and a crosscut in $E$, and suppose $q(c)\notin\overline{\Delta_c}$ for $c\in\Omega$. Fix $p\in\Omega$, a $K$-quasiconformal extension $H_p$ of $h_p$, and a disk-preserving M\"obius map $M_p$ with $M_p(\infty)=H_p^{-1}(q(p))$. For each $0<m\le1$ there is $\delta(K,m)>0$ such that every subarc $J\subset A$ satisfies
\[
 \omega_p(J;\Omega)\ge m\quad\Longrightarrow\quad
 \operatorname{length}(M_p^{-1}\beta(J))\ge\delta(K,m).
\]
The length is the angular length on $\T$.
\end{lemma}
\begin{proof}
Use $F_u$ from Lemma \ref{lem:mainlemma} with $a(c)=q(c)$. The family hypotheses follow from Lemma \ref{lem:family}, and continuity gives $q(c)\notin\Delta_c$ on $\partial\Omega$. Thus $|F_u|\le\eta_K(4)$ there, as in Lemma \ref{lem:mainlemma}. On $J$, where $q(c)=h_c(\beta(c))$, the sharper estimate \eqref{eq:boundary} applies. If $\ell=\operatorname{length}(M_p^{-1}\beta(J))\le1/4$, choose $u$ as in Step~2 of Lemma \ref{lem:mainlemma}. With $\mu=\omega_p(J;\Omega)\ge m$, its harmonic-mean estimate gives
\[
 -\log\eta_K(3)\le\mu\log\eta_K(4\ell)+(1-\mu)\log\eta_K(4)
 \le m\log\eta_K(4\ell)+(1-m)\log\eta_K(4).
\]
This implies the conclusion.
\end{proof}

\subsection{Separation and scale estimates at the level-one root}
Suppose first that $n=1$. Then $c_*=1$, $q(c)=c$, and $\beta=\alpha$.
Let $r(c)$ be the preimage of the critical value $P_c(c)$ under the map $P_c$ except $q(c)=c$.
By calculation,
\begin{equation}\label{eq:level-one-factor}
 P_c(w)-P_c(c)=\frac{\lambda}{3c}(w-c)^2(w-r(c)).
\end{equation}
It follows that $r(c)=(3-c)/2$.

Choose a small quasidisk $\Omega\subset E$ bounded by an arc $A\subset\Gamma$ through $1$ and a crosscut in $E$, with $E\cap B(1,\varepsilon)\subset\Omega$ for some $\varepsilon>0$. For parameters $c\in\Omega$, $c\notin\overline{\Delta_c}$ in the dynamical plane. For $c=\psi(z)\in\mathcal C$ near $1$, the linearizing identity gives
$P_c(c)=h_c(\lambda z)=P_c(h_c(z))$. Since $h_c(z)\ne c$, 
\begin{equation}\label{eq:r-leaf}
 r(\psi(z))=h_{\psi(z)}(z).
\end{equation}
Continuity extends this identity to $z\in\T$ near $1$. Thus $r(c)\in\Delta_c$ for $c\in\mathcal C$ and $r(c)\in\partial\Delta_c$ for $c\in\partial\mathcal C$ near the root. Meanwhile, on $A\setminus\{1\}$ the point $c$ lies on $\partial\Delta_c$ and $r(c)$ lies outside $\overline{\Delta_c}$: otherwise $P_c$ would fail to be injective on $\overline{\Delta_c}$.

\begin{lemma}\label{lem:separation}
There is a constant $b>0$ such that, for parameters sufficiently close to $1$, the following estimates hold.
\begin{enumerate}
    \item $\operatorname{dist}(r(c),\overline{\Delta_c})\ge b|c-1|$ for $c\in A$,
    \item  $\operatorname{dist}(c,\overline{\Delta_c})\ge b|c-1|$ for $c\in\overline{\mathcal C}$.
\end{enumerate}
\end{lemma}

\begin{proof}
For simplicity, let $d=c-1$. Note that $H_c\circ R_\theta\circ H_c^{-1}$ is a $K^2$-quasiconformal map, which agrees with $P_c$ on $\overline{\Delta_c}$. 

Assume (a) failed along a sequence approaching $1$, we can select $y\in\overline{\Delta_c}$ with $y=r(c)+o(d)$, it follows from equation \eqref{eq:level-one-factor},
\[
\begin{aligned}
P_c(y)-P_c(c)&=\frac{\lambda}{3c}(y-c)^2(y-r(c))=O(d^2)\,o(d)=o(d^3),\\
P_c(y)-P_c(1)&=-\frac{\lambda}{6c}d^3+o(d^3).
\end{aligned}
\]
Therefore,
\[
\begin{aligned}
[y,c;1,0]=\frac{y-c}{(y-1)c}
=\frac{-\frac32d+o(d)}{\bigl(-\frac12d+o(d)\bigr)(1+d)}\longrightarrow 3,\\
[P_c(y),P_c(c);P_c(1),0]
=\frac{P_c(y)-P_c(c)}{P_c(y)-P_c(1)}
  \frac{P_c(1)}{P_c(c)}\longrightarrow 0.
\end{aligned}
\]
This contradicts cross-ratio distortion. 

For (b), take $y=c+o(d)$ and use the pair $y,r(c)\in\overline{\Delta_c}$; the limit of the cross ratio $[y,r(c);1,0]$ is $3/2$ and the limit of cross ratio of their images is again zero. Shrinking the neighborhood gives both estimates with one $b>0$.
\end{proof}

\begin{lemma}[Root-scale estimates]\label{lem:root-scale}
Write $p_z=\psi(z)$ for $z\in\T$ near $1$. There are constants with the following properties:
\begin{enumerate}[label=(\roman*),leftmargin=*]
\item $\operatorname{dist}(p_z,\Gamma)\ge\epsilon |p_z-1|$ for $z\ne1$.
\item $|z-1|\le|w-1|$ implies $|p_z-1|\le C |p_w-1|$.
\item For each $a_0>0$, if $|z-w|\ge a_0\max\{|z-1|,|w-1|\}$, then
\[
 |p_z-p_w|\ge b_{a_0}\max\{|p_z-1|,|p_w-1|\}.
\]
\end{enumerate}
\end{lemma}
\begin{proof}
On $A$ the critical-angle coordinate of $a(c)=c$ is $\alpha(c)$, with $\alpha(1)=1$. It is quasisymmetric by Theorem~\ref{thm:main} and the quasicircularity of $\Gamma$. Fix a Riemann map $g:\D\to\Omega$ with $g(1)=1$. Its boundary extension is quasisymmetric, and so is $\tau=\alpha\circ g$ on a boundary neighborhood of $1$. For $z\in\T$, $p_z\in\partial\mathcal{C}\cap\Omega$. 

The normalized extensions $H_p$ and their inverses have uniform local quasisymmetry near $1$. For simplicity, we write $p=p_z$ and $\xi_p=H_p^{-1}(p)$. By definition, $H_p(1)=1$, $H_p(\xi_p)=p$. By equation \eqref{eq:r-leaf}, $H_p(z)=r(p)$. By calculation, $\frac{|p-1|}{|r(p)-1|}=2$, in the dynamical plane. By the quasisymmetry of $H_p^{-1}$, we obtain
\begin{equation}\label{eq:xi-scale}
 |\xi_p-1|\asymp |z-1|.
\end{equation}

Since $H_p(\xi_p/|\xi_p|)\in \partial\Delta_p$, by Lemma \ref{lem:separation} (b), $|p-H_p(\xi_p/|\xi_p|)|\ge b|p-1|$.
Again by the quasisymmetry of $H_p^{-1}$, $|\xi_p-1|\le \eta(1/b)|\xi_p-\xi_p/|\xi_p||$.

By the triangle inequality, we obtain
\begin{equation}\label{eq:xi-scale-1}
|\xi_p|-1\asymp|\xi_p-1|.
\end{equation}

Let $M_p$ preserve $\D$ and send $\infty$ to $\xi_p$, thus $M_p^{-1}(v)=e^{i\gamma_p} \frac{v-1/\overline{\xi_p}}{1-v/\xi_p}$. On $\T$,
\begin{equation}\label{eq:mobius-density}
|(M_p^{-1})'(v)|=\frac{|\xi_p|^2-1}{|v-\xi_p|^2},
\end{equation}
and by equation \eqref{eq:xi-scale-1} and triangle inequality,
\begin{equation}\label{eq:mobius-density-1}
 |v-\xi_p|\asymp |v-1|+|\xi_p-1|.
\end{equation}

For (i), suppose along a sequence that $\operatorname{dist}(p,\Gamma)/|p-1|\to0$. Quasisymmetry of the extension of $g$ implies that $u=g^{-1}(p)\in\D$ approaches $1$ tangentially. Take a boundary arc $J_0$ of angular radius $1-|u|$ centered at $u/|u|$. Its harmonic measure at $u$ is bounded below, whereas
 $\frac{\operatorname{diam}(J_0)}{\operatorname{dist}(J_0,1)}\longrightarrow0$, since $u$ approaches $1$ tangentially. Local quasisymmetry gives the same ratio limit for $\tau(J_0)$. Equations~\eqref{eq:mobius-density} and \eqref{eq:mobius-density-1} now give
\[
\begin{aligned}
\operatorname{length}(M_p^{-1}(\tau(J_0)))
&\le \operatorname{length}(\tau(J_0))\sup_{v\in \tau(J_0)}|(M_p^{-1})'(v)| \\
&\le C\operatorname{length}(\tau(J_0))\,\frac{|\xi_p-1|}{(\operatorname{dist}(\tau(J_0),1)+|\xi_p-1|)^2} \\
&\le C'\frac{\operatorname{length}(\tau(J_0))}{\operatorname{dist}(\tau(J_0),1)}
\longrightarrow 0.
\end{aligned}
\]
This contradicts Lemma~\ref{lem:local-mass} applied to $a(c)=c$, and proves (i).

For (ii), let $s_z=|g^{-1}(p_z)-1|$ and $J_z=\{e^{it}:s_z\le t\le2s_z\}$. 

By (i), $1-|g^{-1}(p_z)|\asymp s_z$,
so $\omega_{g^{-1}(p_z)}(J_z;\D)\ge m_0>0$. Quasisymmetry of $\tau$ implies for $v\in J_z$,
\[
 \operatorname{length}\tau(J_z)\asymp|\tau(v)-1|
 \asymp|\tau(e^{is_z})-1|.
\]
Hence Lemma~\ref{lem:local-mass} and \eqref{eq:xi-scale}--\eqref{eq:mobius-density-1} imply
\[
 0<\delta_0\le\operatorname{length}(M_{p_z}^{-1}\tau(J_z))
 \le C\frac{|z-1|\,|\tau(e^{is_z})-1|}
 {\bigl(|z-1|+|\tau(e^{is_z})-1|\bigr)^2}.
\]
Since $ab/(a+b)^2\le\min\{a/b,b/a\}$ for $a,b>0$, this yields
\begin{equation}\label{eq:angle-scale}
 |z-1|\asymp|\tau(e^{is_z})-1|.
\end{equation}
If $|z-1|\le|w-1|$, quasisymmetry of $\tau^{-1}$ now gives $s_z\le C s_w$. Since $|e^{is_z}-1|\asymp s_z=|g^{-1}(p_z)-1|$, quasisymmetry of $g$ yields
\[
 |p_z-1|\asymp|g(e^{is_z})-g(1)|
 \le C|g(e^{is_w})-g(1)|\asymp|p_w-1|,
\]
proving (ii).

For (iii), fix $a_0>0$. Assume that there is a sequence of pairs $z,w\ne1$ satisfying the assumption in (iii) but $|p_z-p_w|/\max\{|p_z-1|,|p_w-1|\}\to0$. The triangle inequality implies $|p_z-1|/|p_w-1|\to1$. Since the crosscut part of $\partial\Omega$ stays away from $1$, (i) implies that $B(p_w,\kappa|p_w-1|)\subset\Omega$ for a fixed $\kappa>0$. Writing $d_U$ for hyperbolic distance in a domain $U$, the disk estimate therefore gives
\[
 d_\Omega(p_z,p_w)\le C\frac{|p_z-p_w|}{|p_w-1|}\longrightarrow0.
\]
For each fixed $z\ne1$, let $G_z(c)=(h_c(z)-1)/(c-1)$. By Lemma~\ref{lem:family}, this function is holomorphic in $\Omega$. It omits $0$ because $h_c$ is injective on $\overline\D$ and $h_c(1)=1$; it omits $1$ because $c\notin\overline{\Delta_c}$. Thus $G_z:\Omega\to\C\setminus\{0,1\}$. By equation \eqref{eq:r-leaf}, $G_z(p_z)=-1/2$, and $h_{p_w}(w)-1=-(p_w-1)/2$. The Schwarz–Pick lemma yields
\[
 d_{\C\setminus\{0,1\}}\bigl(G_z(p_w),-1/2\bigr)
 \le d_\Omega(p_w,p_z)\longrightarrow0.
\]
Since the base point $-1/2$ is fixed, this implies
\[
 \frac{h_{p_w}(z)-h_{p_w}(w)}{p_w-1}
 =G_z(p_w)+1/2\longrightarrow0.
\]
On the other hand, $|z-w|\ge a_0|w-1|$ by assumption. Applying the uniform local quasisymmetry of $h_{p_w}$ with base point $w$, and using $|h_{p_w}(w)-1|=|p_w-1|/2$, we obtain
\[
 \frac{|h_{p_w}(z)-h_{p_w}(w)|}{|p_w-1|}
 \ge\frac{1}{2\eta(1/a_0)}>0,
\]
where $\eta$ is independent of $p_w$. This contradiction proves (iii).
\end{proof}

\subsection{Small boundary arcs at the level-one root}
Continue with the fixed level-one component $\mathcal{C}$ and the domain $\Omega$ above. This is the unique level-one component in $E$ \cite[Theorem 5.2]{WY23}, which simplifies the small-arc argument below.

\begin{lemma}[Small arcs at the level-one root]\label{lem:small-arcs}
For the fixed component $\mathcal{C}$ with $n=1$ and $c_*=1$, there exist $a>0$ and $C_0<\infty$ such that every arc $J\subset\T$ sufficiently close to $1$, with endpoints $z,w$ satisfying 
\[
 \sup_{v\in J}|v-z|\le a|z-1|,
\]
also satisfying $\operatorname{diam}\psi(J)\le C_0|\psi(z)-\psi(w)|$.
\end{lemma}
\begin{proof}
Recall $r(c)=(3-c)/2$. For $z\in\overline\D\setminus\{1\}$ near $1$, let
\begin{equation}\label{eq:B-pencil}
 B_z(c)=[c,1;h_c(z),r(c)].
\end{equation}
This is holomorphic in $c\in\Omega$. For $c\in A$, by Lemma \ref{lem:separation} (a), $|h_c(z)-r(c)|\ge \operatorname{dist}(r(c),\overline{\Delta_c})\ge b|c-1|$. Therefore,
\[
|B_z(c)|
=\frac{2|r(c)-h_c(z)|}{|c-h_c(z)|}\ge
\frac{2|r(c)-h_c(z)|}
{|r(c)-h_c(z)|+\frac32|c-1|}\ge\frac{4b}{2b+3}>0.
\]
For $c$ on the crosscut, $B_z\to B_1=1$ uniformly as $z\to1$. As $c\to1$, for fixed $z\ne1$, since $h_1(z)\not=h_1(1)=1$, $B_z\to \frac{2(h_1(z)-1)}{1-h_1(z)}=-2$.  Hence there is a constant $b_0>0$ such that for any $c\in\partial\Omega$,
\begin{equation}\label{eq:B-boundary}
 |B_z(c)|\ge b_0>0\qquad.
\end{equation}

We consider the zeros of $B_z(c)=\frac{2(h_c(z)-r(c))}{c-h_c(z)}$ on $\Omega$. Since $c\notin\overline{\Delta_c}$ for $c\in \Omega$, it suffices to consider the zero of $h_c(z)-r(c)$ on $\Omega$. For $v\in\D$, a zero $c_0$ of $r(c)-h_c(v)$ satisfies
\[
P_{c_0}(c_0)=P_{c_0}(r(c_0))=P_{c_0}(h_{c_0}(v))=h_{c_0}(\lambda v)\in\Delta_{c_0}.
\]
It therefore lies in the unique level-one component, with $\Phi_{\mathcal C}(c_0)=\lambda v$. This parameter is unique because  $\Phi_{\mathcal C}$ is conformal on $\mathcal{C}$. Meanwhile, for any $c\in \mathcal C$,
\[
P_c(r(c))-P_c(h_c(v))
=h_c(\Phi_{\mathcal C}(c))-h_c(\lambda v).
\]
Take the derivative of both sides of the equation with respect to $c$ at $c=c_0$. Since $P_c'(r(c))=\frac{3\lambda(c-1)^2}{4c}\not=0$ and $h_{c_0}$, $\Phi_{\mathcal{C}}$ are conformal, for $c\in \Omega$,
\[
\left.\frac{d}{dc}\bigl(r(c)-h_c(v)\bigr)\right|_{c=c_0}\ne0.
\]
This implies the zero $c_0$ is simple. Take $z\in\T\setminus\{1\}$, and let $v=tz$ for $0<t<1$. For $t\to1$, 
\[
r(c)-h_c(tz)\longrightarrow r(c)-h_c(z)
\]
uniformly on every compact subset of \(\Omega\).
The limit function is not identically zero, since as \(c\to1\),
\[
r(c)-h_c(z)\longrightarrow 1-h_1(z)\neq 0.
\]
By Rouché's theorem, for $z\in \T\setminus\{1\}$, $r(c)-h_c(z)$ has at most one zero, counted with multiplicity. Since $r(\psi(z))=h_{\psi(z)}(z)$, and for $z\not=1$ near $1$, $p=\psi(z)\in \Omega$, thus for $z\in\T\setminus\{1\}$ near $1$, $B_z$ has precisely one zero in $\Omega$ which is $p=\psi(z)$, and this zero is simple.

Choose $g_p:\D\to\Omega$ conformal with $g_p(0)=p$, and set $F_z=B_z\circ g_p$. The function $u/F_z(u)$ is holomorphic on $\D$ and continuous on $\overline\D$. The maximum principle and \eqref{eq:B-boundary} imply
\begin{equation}\label{eq:level-one-lower}
 |F_z(u)|\ge b_0|u|.
\end{equation}

Over $u\in\D$, and $h_p(v)=w\in\overline{\Delta_p}\cup\{p,\infty\}$ define
\[
\mathcal H_u(w)=
\begin{cases}
h_{g_p(u)}\!\left(h_p^{-1}(w)\right),
& w\in\overline{\Delta_p},\\[2mm]
g_p(u),&w=p,\\
\infty,&w=\infty.
\end{cases}
\]
By Lemma \ref{lem:family} and $g_p(u)\notin\overline{\Delta_{g_p(u)}}$,
$\mathcal H_u(w)$ is a holomorphic motion.
Extend it to a sphere motion Hu using Słodkowski’s theorem \cite{Slo91}; its dilatation is at most 3 on $\D_{1/2}$. 
Let $T_z(w)=[p,1;r(p),w]$ be a M\"{o}bius transformation such that $T_z(r(p))=0$, $T_z(1)=\infty$ and $T_z(p)=1$.
Define
\begin{equation}
Q_z(u)=T_z\bigl(\mathcal H_u^{-1}(r(g_p(u)))\bigr).
\end{equation}
The holonomy theorem \cite[Lemma 17.9]{Lyu22} implies that $Q_z$ is a $3$-quasiregular map on $\D_{1/2}$. It satisfies $Q_z(0)=0$, has no poles since $r(c)\ne1$ in $\Omega$, and is nonconstant by \eqref{eq:B-boundary}. Cross-ratio distortion for $\mathcal H_u$ gives
\[
 |F_z(u)|\le\eta_3(|Q_z(u)|),\quad \text{for} \quad u\in \D_{1/2}.
\]
Therefore, 
\[
 |Q_z(u)|\ge\sigma:=\eta_3^{-1}(b_0/4)>0\quad \text{for} \quad |u|=1/4.
\]

Fix $0<R<\sigma$. The component $V$ of $Q_z^{-1}(\D_R)$ containing $0$ is compactly contained in $\D_{1/4}$, and $Q_z:V\to\D_R$ is proper. We claim that $Q_z:V\to\D_R$ is of degree one. Take $y\in\D_R\cap T_z(\Delta_p)$, such that
$v=h_p^{-1}(T_z^{-1}(y))\in\mathbb D$. $\D_R\cap T_z(\Delta_p)\not=\emptyset$ since 
$r(p)=h_p(z)\in\partial\Delta_p$ and $T_z(r(p))=0$, i.e., $0$ is a boundary point of $T_z(\Delta_p)$. Let \(c=g_p(u)\), it follows that
\[
\begin{aligned}
Q_z(u)=y
&\iff \mathcal H_u^{-1}(r(c))=T_z^{-1}(y)=h_p(v)\\
&\iff r(c)=\mathcal H_u(h_p(v))\iff r(c)=h_c(v).
\end{aligned}
\]
Therefore, the preimages of $y$ correspond to zeros of $r(c)-h_c(v)$ with $v\in\D$. There is at most one such zero. On the capture component, since
$\mathcal H_u^{-1}(r(c))
=h_p\!\left(\lambda^{-1}\Phi_{\mathcal C}(c)\right)$, 
\begin{equation}\label{eq:Q-on-capture}
 Q_z(u)=T_z\bigl(h_p(\lambda^{-1}\Phi_{\mathcal C}(g_p(u)))\bigr).
\end{equation}
All the maps on the right-hand side are conformal here, so the local degree is one. Let $U$ be the component of $Q_z^{-1}(\D_{R/2})$ containing $0$.

Now assume the condition $|v-z|\le a|z-1|$ with $v\in J$. Since $p-r(p)=3(p-1)/2$, uniform local quasisymmetry of $h_p$ near $1$ gives
\begin{equation}\label{eq:small-T}
 |T_z(h_p(v))|=\frac23\frac{|h_p(z)-h_p(v)|}{|1-h_p(v)|}
 \le\frac23\eta\!\left(\frac a{1-a}\right)<R/2\quad \text{for}\quad v\in J ,
\end{equation}
by taking $a$ sufficiently small. 
By definition, 
\[
Q_z(g_p^{-1}(\psi(v)))=T_z(\mathcal H_u^{-1} (r(\psi(v))))=T_z(\mathcal H_u^{-1}(h_{\psi(v)}(v))) = T_z(h_p(v)).
\]
By equation \eqref{eq:small-T}, 
\[
Q_z\!\left(g_p^{-1}(\psi(J))\right)
=T_z(h_p(J))
\subset\mathbb D_{R/2}.
\]

The curve $\partial\Delta_p$ is a uniform $K$-quasicircle and $T_z(1)=\infty$. Thus $T_z(\partial\Delta_p\setminus\{1\})$ is a uniform quasiline, that is, the image of $\mathbb R$ under a global $K$-quasiconformal map fixing $\infty$. Global quasisymmetry gives a uniform turning bound $C(K)$ for its finite subarc $T_z(h_p(J))$.

Factor $Q_z\circ g_p^{-1}=f\circ L$ on $g_p(\D_{1/2})$, where $L$ is a global $3$-quasiconformal map fixing $\infty$, such that $f:L(g_p(V))\to\D_R$ is conformal. Since $f(L(\psi(J)))=T_z(h_p(J))$ has turning bounded by $C(K)$ in $\D_{R/2}$, the turning-distortion lemma \cite[Lemma 6.1(2)]{QWY12} implies
\[
\frac{\operatorname{diam}L(\psi(J))}
 {|L(\psi(z))-L(\psi(w))|}
\le
D(d,m)
\frac{\operatorname{diam}T_z(h_p(J))}
 {|T_z(h_p(z))-T_z(h_p(w))|}\le C(K)D(d,m),
\]
where $d=1$ and $m=\operatorname{mod}
\bigl(\mathbb D_R\setminus\overline{\mathbb D_{R/2}}\bigr)=\frac{\log2}{2\pi}.$ Global quasisymmetry of $L^{-1}$ gives
\[
\operatorname{diam}\psi(J)
\le C_1|\psi(z)-\psi(w)|.
\]
\end{proof}

\begin{proof}[Proof of Theorem \ref{thm:capture} at the level-one root]
We now combine the two scales. Let $J\subset\T$ be a sufficiently short arc near $1$, with endpoints $z,w$. If $1$ is an endpoint, for all $v\in J$, $|v-1|\le|w-1|$, and $|\psi(v)-1|\le C|\psi(w)-1|$ by Lemma~\ref{lem:root-scale}(ii). Therefore, for any pair $v_1,v_2\in J$,
\[
\operatorname{diam}\psi(J)=\sup|\psi(v_1)-\psi(v_2)|\le|\psi(v_1)-1|+|\psi(v_2)-1|\le 2C|\psi(w)-1|.
\]
Let $s=\max\{|z-1|,|w-1|\}$. If $1$ is not an endpoint, without loss of generality we assume the endpoint $z$ satisfying $|z-1|=s$. If $\sup_{v\in J}|v-z|\le as$, Lemma~\ref{lem:small-arcs} applies. In the other case, $|z-w|=\sup_{v\in J}|v-z|\ge as$. Since $J$ lies in a small circular neighborhood of $1$, every $v\in J$ satisfies $|v-1|\le s$. By Lemma~\ref{lem:root-scale}(ii),
\[
\operatorname{diam}\psi(J)\le 2C|p_z-1|\le 2C\max\{|p_z-1|,|p_w-1|\}.
\]
Therefore Lemma~\ref{lem:root-scale}(iii) implies
\begin{equation}\label{eq:root-turning}
 \operatorname{diam}\psi(J)\le C'\max\{|p_z-1|,|p_w-1|\}
 \le C''|\psi(z)-\psi(w)|.
\end{equation}
Thus the level-one capture boundary is locally a quasiarc at $1$ by Theorem \ref{thm:bounded-turning}.
\end{proof}

\subsection{Uniqueness and boundary geometry at higher-level roots}
Now assume that the fixed attached component has level $n\ge2$. We first prove uniqueness at its root, then verify the hypotheses of the level-one argument.

\begin{lemma}[Uniqueness at an attachment point]\label{lem:unique-root}
The fixed root $c_*$ belongs to the closure of exactly one capture component in $E$, namely $\mathcal C$. Moreover, $q'(c_*)\ne0$.
\end{lemma}
\begin{proof}
By \cite{WY23}, for $k\ge1$ and $\zeta\in\D$, the equation
\[
 P_c^k(c)=h_c(\zeta)
\]
has exactly $(3^k-1)/2$ simple zeros in $E$: one in each component $\mathcal B$ of level $m\le k$, given by $\Phi_{\mathcal B}(c)=\lambda^{m-k}\zeta$. Each zero has a limit as $\zeta$ tends radially to $\T$, since $\Phi_{\mathcal B}^{-1}$ extends homeomorphically to $\overline\D$. Rouch\'e's theorem on small disks compactly contained in \(E\) counts the interior limits with multiplicity.

The polynomial $q(c)=P_c^{n-1}(c)$ has degree $D=(3^{n-1}+1)/2$ by calculation. Apply the preceding observation to $q(c)=h_c(s)$ as $s\nearrow1$. None of its $D-1$ zeros can limit on $\Gamma$: a limit in the closure of a level-$m$ component, $m<n$, would satisfy both $P_c^{m-1}(c)=1$ and $q(c)=1$, making $1$ periodic. Thus $q(c)-1$ has at least $D-1$ zeros in $E$, counted with multiplicity. Together with $q(c_*)=1$, this exhausts its degree. Consequently $c_*$ is simple and is the only zero of $q(c)-1$ outside $E$.

Set $A(c)=c(P_c^n(c)-P_c(1))$. Direct calculation gives
\begin{equation}\label{eq:root-count-factor}
 A(c)=\frac{\lambda}{3}(q(c)-1)^2
 \left(q(c)-\frac{3c-1}{2}\right).
\end{equation}
It has degree $3D$ and $A(0)=\lambda/6\ne0$. On $\Gamma$, injectivity of $P_c$ on $\partial\Delta_c$ implies that $A(c)=0$ exactly when $q(c)=1$. Hence its only zero on $\Gamma$ is $c_*$, of multiplicity two; here $c_*\ne1$, since $n\ge2$.

We claim that there is exactly one simple zero of $A$ in $I$. For $c\in I\setminus\{0\}$, inversion symmetry gives
\[
 P_c^n(c)-P_c(1)
 =-c\bigl(P_{1/c}(1/c)-h_{1/c}(\lambda^n)\bigr).
\]
For $0<s<1$, the equation $P_a(a)=h_a(s\lambda^n)$ has one simple zero in $E$, in the unique level-one component. Its limit stays in $E$, since a limit on \(\Gamma\) would have boundary coordinate \(\lambda\), whereas its coordinate tends to \(\lambda_n \neq \lambda\). Rouch\'e's theorem therefore gives exactly one simple zero of $P_a(a)-h_a(\lambda^n)$ in $E$. Thus $A$ has $3D-3$ zeros in $E$, counted with multiplicity.

Finally, $P_c^n(c)=h_c(s\lambda)$ has $(3^n-1)/2=3D-2$ simple zeros in $E$. As $s\nearrow1$, exactly $3D-3$ of them converge to the zeros of $A$ in $E$. The remaining one tends to $c_*$. Every level-$n$ component attached to $c_*$ contributes such a zero, so there is exactly one. A component of another level cannot share this root, since the two root relations would again make $1$ periodic.
\end{proof}

\emph{The local partner.}
It is clear that $c_*\ne1$. The expansion
\[
 P_c(1+u)=P_c(1)+\frac{\lambda(1-c)}{2c}u^2+\frac{\lambda}{3c}u^3
\]
has a nonzero quadratic coefficient at $c_*$. For a point $q(c)$ near the critical point $1$, one can holomorphically choose another point $r(c)$ near $1$ such that $q(c)$ and $r(c)$ have the same image under $P_c$ with
\begin{equation}\label{eq:partner}
 P_c(r(c))=P_c(q(c)),\quad \text{and} \quad
 r(c)-1=-(q(c)-1)+O((q(c)-1)^2).
\end{equation}
In particular $t(c)=(r(c)-1)/(q(c)-1)$ extends holomorphically across $c_*$ with $t(c_*)=-1$.

Choose a quasidisk $\Omega\subset E$ bounded by an arc $A\subset\Gamma$ through $c_*$ and a crosscut, containing $E\cap B(c_*,\varepsilon)$, and small enough that
\begin{equation}\label{eq:high-exterior}
\begin{gathered}
    q(c) \notin \overline{\Delta_c} \quad \text{for} \quad c \in \Omega,\quad \text{and} \\
    q(c) - 1 \neq 0, \quad r(c) - 1 \neq 0, \quad q(c) - r(c) \neq 0
    \quad \text{for} \quad c \in \overline{\Omega} \setminus \{c_*\}.
\end{gathered}
\end{equation}
Every $c\in E$ with $q(c)\in\overline{\Delta_c}$ lies in the closure of a component of level less than $n$. For $q(c)\in\Delta_c$ this is immediate; for a boundary value $q(c)=h_c(\zeta)$, apply Rouch\'e's theorem to $q-h_c(s\zeta)$ as $s\nearrow1$. Its limit is not identically zero, since $q(c)\to\infty$ at infinity while $h_c(\zeta)$ stays bounded. There are finitely many component closures of levels at most $n$, and Lemma~\ref{lem:unique-root} shows that only $\overline{\mathcal C}$ contains $c_*$. Shrinking $\Omega$ to avoid the others gives the first condition and ensures that every level-$n$ component meeting $\Omega$ is $\mathcal C$. The remaining nonvanishing conditions follow from \eqref{eq:partner} and the simple zero of $q-1$ at $c_*$.

For $c=\psi(z)$ near $c_*$, the two local solutions of $P_c(w)=h_c(\lambda z)$ are $q(c)$ and $r(c)$. By equation \eqref{eq:high-exterior},
\begin{equation}\label{eq:high-leaf}
 r(\psi(z))=h_{\psi(z)}(z)\quad\text{for} \quad z\in\overline\D\text{ near }1.
\end{equation}
On $A$, $q(c)\in\partial\Delta_c$ and $r(c)\notin\overline{\Delta_c}$ when $c\not=c_*$. The proof of Lemma~\ref{lem:separation} now yields
\begin{align}
 \operatorname{dist}(r(c),\overline{\Delta_c})&\ge b|q(c)-1|\quad  \text{for}\quad c\in A,\label{eq:high-gap-gamma}\\
 \operatorname{dist}(q(c),\overline{\Delta_c})&\ge b|q(c)-1|\quad \text{for}\quad c\in\overline{\mathcal C}\text{ near }c_*.\label{eq:high-gap-capture}
\end{align}
Indeed, let $d=q(c)-1$. The quadratic expansion gives $|P_c(q(c))-P_c(1)|\asymp|d|^2$. For $y$ within $\varepsilon_0|d|$ of $r(c)$ or $q(c)$, respectively, the original cross ratios are uniformly bounded away from zero, whereas the image cross ratios are bounded by $C|y-r(c)|/|d|$ or $C|y-q(c)|/|d|$. Cross-ratio distortion gives both estimates exactly as in Lemma~\ref{lem:separation}.

\emph{The root-scale estimates.}
The estimates in Lemma~\ref{lem:root-scale}(i)--(iii) also hold at $c_*$, with $p_z=\psi(z)$ and the parameter distances $|p_z-1|$, $|p_w-1|$ replaced by $|p_z-c_*|$, $|p_w-c_*|$. To see this, use the following functions in its proof:
\[
 \beta(c)=\lambda^{n-1}\alpha(c),\quad
 \xi_p=H_p^{-1}(q(p)),\quad
 G_z(c)=\frac{h_c(z)-1}{q(c)-1},\quad \text{and} \quad g(1)=c_*.
\]
The coordinate $\beta$ on $A$ is quasisymmetric and satisfies $\beta(c_*)=1$. Equations \eqref{eq:partner}, \eqref{eq:high-leaf}, and \eqref{eq:high-gap-capture} give \eqref{eq:xi-scale}, and Lemma~\ref{lem:local-mass} applies to $a=q$. Thus the proofs of (i) and (ii) are unchanged. In (iii), $G_z(p_z)=t(p_z)\to-1$ replaces the constant $-1/2$. Continuity and Schwarz--Pick Lemma give
\[
 \frac{h_{p_w}(z)-h_{p_w}(w)}{q(p_w)-1}
 =G_z(p_w)-t(p_w)\longrightarrow0
\]
under the same contrary assumption, contradicting uniform local quasisymmetry as before.

\emph{Small arcs at the level-n root.} Repeat the proof of Lemma~\ref{lem:small-arcs} with $q,r,\Omega$ constructed above and
\begin{equation*}\label{eq:high-B-pencil}
 B_z(c)=[q(c),1;h_c(z),r(c)]\quad\text{and}\quad
 T_z(w)=[q(p),1;r(p),w].
\end{equation*}
In the holomorphic motion, replace the moving marked point $p\mapsto g_p(u)$ with $q(p)\mapsto q(g_p(u))$. The boundary lower bound \eqref{eq:B-boundary} follows from \eqref{eq:high-gap-gamma}, with root limit $-1$ replacing $-2$.

\begin{proof}[Completion of the proof of Theorem~\ref{thm:capture}]
For the fixed higher-level component, use the root-scale estimates and the small-arc estimates at the level-n root above. Hence $\partial\mathcal C$ is locally a quasiarc at $c_*$, as already proved at the level-one root. Wang and Yang's theorem shows the quasiarc property at every other boundary point. Compactness gives the bounded-turning condition, so Theorem~\ref{thm:bounded-turning} makes $\partial\mathcal C$ a quasicircle. This completes the proof.
\end{proof}

\subsection*{Acknowledgments} This work was supported by the National Natural Science Foundation of China (Grant No.~12526526).

\bibliographystyle{amsalpha}
\bibliography{Ref1}

\end{document}